\documentclass[11pt]{article}

\usepackage[letterpaper,margin = 1in]{geometry}

\usepackage[english]{babel}

\usepackage{amssymb,mathrsfs,amsmath,amsfonts,verbatim, amsthm}
\usepackage{graphicx}
\usepackage{mathtools}
\usepackage{bm,bbm}
\usepackage[dvipsnames]{xcolor}
\usepackage{scalerel}
\usepackage{showtags}
\usepackage{enumerate}
\usepackage{microtype}
\usepackage{stackengine}
\usepackage{cmap}
\usepackage[colorlinks=true,allcolors=blue]{hyperref}
\usepackage[normalem]{ulem}

\usepackage{orcidlink}

\newtheorem{theorem}	{Theorem}[section]

\newtheorem{proposition}[theorem]{Proposition}
\newtheorem{lemma}		[theorem]{Lemma}

\newtheorem{conjecture}	[theorem]{Conjecture}

\theoremstyle{definition}

\newcommand{\C}{{\mathbb{C}}}

\newcommand{\D}{{\mathbb{D}}}

\newcommand{\bbT}{\mathbb{T}}

\title{Optimal Polynomial Approximants vanishing arbitrarily close to the origin}

\author{Athanasios Beslikas \orcidlink{0009-0006-8098-6843}, Christopher Felder \orcidlink{0000-0002-0736-0987}, and Douglas T. Pfeffer \orcidlink{0000-0001-8287-2686}}
\date{}

\begin{document}
	
	\maketitle
	
	\begin{abstract}
		This work concerns a sequence of stable polynomials on the bidisk whose affine optimal polynomial approximants in the Hardy space of the bidisk have zero varieties that come arbitrarily close to the origin.
	\end{abstract}

	\renewcommand{\thefootnote}{\fnsymbol{footnote}} 
	\renewcommand{\thefootnote}{\arabic{footnote}} 
	
	\renewcommand{\thefootnote}{\fnsymbol{footnote}} 
	\renewcommand{\thefootnote}{\arabic{footnote}}

	\newcommand{\Addresses}{{
			\bigskip
			\footnotesize
			
			\setlength{\parindent}{0cm}
			\noindent Athanasios Beslikas \,\orcidlink{0009-0006-8098-6843}\par\nopagebreak  \textsc{Doctoral school of exact and natural studies, Institute of Mathematics, Faculty of Mathematics and Computer Science, Jagiellonian University, Lojasiewicza 6, 30-348, Krakow}\par\nopagebreak
			e-mail: \texttt{athanasios.beslikas@doctoral.uj.edu.pl}
			
			\medskip
			
			Christopher Felder \,\orcidlink{0000-0002-0736-0987}\par\nopagebreak  \textsc{Department of Mathematics \& Statistics, University of South Florida, 4202 E. Fowler Avenue, Tampa, FL 33620
			}\par\nopagebreak
			e-mail: \texttt{felderc@usf.edu}
			
			\medskip
			
			Douglas T. Pfeffer\,\orcidlink{0000-0001-8287-2686}\par\nopagebreak\textsc{Department of Mathematics, University of Tampa, 401 W. Kennedy Blvd. Tampa, FL 33606}\par\nopagebreak
			e-mail: \texttt{dpfeffer@ut.edu}
	}}

	\section{Introduction and Background}
	\label{sec:introduction}
	
	Let $\mathbb{D}=\{z\in\mathbb{C}:|z|<1\}$ denote the open unit disk and let $\D^d$ denote the $d$-fold Cartesian product of the unit disk. This work is set in the Hardy space on these domains, 
	\[
	H^2(\D^d) := \left\{ f \in \operatorname{Hol}(\D^d) :   \sup_{0\le r < 1} \int_{\bbT^d} |f(r\zeta)|^2 \, d\mu(\zeta) < \infty \right\},
	\]
	where $\operatorname{Hol}(\D^d)$ is the set of holomorphic functions on the $d$-dimensional polydisk and $d\mu$ is normalized Lebesgue measure on the $d$-dimensional torus $\bbT^d$. It is well known that Hardy space functions have radial limits existing almost everywhere on $\bbT^d$, and
	$H^2(\mathbb{D}^d)$ may be identified with the closed subspace of $L^2(\mathbb{T}^d, d\mu)$ consisting of elements with vanishing negative Fourier coefficients. 
	We will primarily work in the case of the bidisk, where, for $f, g \in \operatorname{Hol}(\D^2)$ given by
	\[
	f(z,w)=\sum_{j,k\geq 0}a_{j,k}z^jw^k, \quad g(z,w)=\sum_{j,k\geq 0}b_{j,k}z^jw^k,
	\]
	we write the Hardy-space inner product interchangeably as
	\[
	\langle f, g \rangle = \sum_{j,k\geq 0}a_{j,k} \overline{b_{j,k}}
	= \int_{\bbT^2} f(\zeta) \overline{g(\zeta)} \, d\mu(\zeta).
	\]
	For general theory of Hardy spaces on polydisks see, e.g., \cite{RudinPolydiscs}. 
	
	Given a nonzero function $f\in H^2(\mathbb{D}^2)$ and a finite-dimensional space of analytic polynomials $\mathcal{P}$, the \textit{optimal polynomial approximant}, or OPA, to $1/f$ from $\mathcal{P}$ is the unique polynomial $p\in\mathcal{P}$ satisfying
	\[
	\|pf-1\|_2=\inf_{q\in\mathcal{P}}\|qf-1\|_2.
	\]
	When
	\[
	\mathcal{P}=\operatorname{span}\{1,z,w\},
	\]
	we call $p$ the \emph{affine optimal polynomial approximant} to $1/f$. Note that such a polynomial exists, is unique, and, as a multiple of $f$, is given by the orthogonal projection of the unit constant onto the finite-dimensional space $\mathcal{P}f$. When $d=1$, OPAs are defined analagously by choosing $\mathcal{P}$ to be a space of one-variable polynomials of a fixed degree.

	Optimal polynomial approximants arise naturally in the study of cyclic vectors for shift operators. A function $f\in H^2(\mathbb{D}^d)$ is called cyclic if the polynomial multiples of $f$ are dense in $H^2(\mathbb{D}^d)$. Since the constant function $1$ is cyclic, this is equivalent to the existence of a sequence of polynomials $p_n$ such that
	\[
	\|p_nf-1\|_2\longrightarrow 0.
	\]
	An OPA determines the polynomial multiple of $f$ that best approximates $1$ in the chosen space of polynomial multiples of $f$. Optimal approximants are connected to invariant subspaces, orthogonal polynomials, and reproducing kernels; see, for example, \cite{BeneteauCondoriLiawSecoSolaI,BeneteauCondoriLiawSecoSolaII,BeneteauKhavinsonLiawSecoSola,BeneteauKhavinsonLiawSecoSimanek,ChengFelder,FelderGeneral,SargentSola}.
	
	The same extremal problem appeared earlier in engineering under the name \emph{least-squares polynomial inversion}. Its development was motivated by prediction theory, stationary stochastic processes, and the design of recursive digital filters. Early antecedents include the work of Robinson on stationary stochastic processes \cite{Robinson}, while Chui studied least-squares inverses in one and several variables and their approximation-theoretic properties \cite{ChuiDouble}. Chui and Chan subsequently used the zero-free behavior of one-variable least-squares inverses in recursive filter design \cite{ChuiChan}. Izumino later reformulated the one-variable problem in operator-theoretic terms \cite{Izumino}. For a modern account of the relation between optimal polynomial approximants and filter design, see \cite{BeneteauCentner,SargentSola}.
	
	\subsection*{The History of the Shanks Conjecture}

	A central question in this theory concerns the location of the zeros of optimal polynomial approximants. In one variable, OPAs in the Hardy space $H^2(\mathbb{D})$ are zero-free in the closed unit disk. This is closely related to the theory of orthogonal polynomials and reproducing kernels. Modern extensions examine the location and accumulation of zeros of OPAs in general weighted Hilbert spaces of analytic functions \cite{BeneteauKhavinsonLiawSecoSola,BeneteauKhavinsonLiawSecoSimanek}.
	
	The corresponding problem in two variables arose from questions about the stability of two-dimensional recursive filters. In 1972, Shanks, Treitel, and Justice proposed that the least-squares inverse polynomials associated with a two-variable polynomial should be zero-free in the closed bidisk \cite{ShanksTreitelJustice}. In the language now used in function theory, their conjecture predicted that if $f$ is a polynomial, then the OPAs to $1/f$ in $H^2(\mathbb{D}^2)$ do not vanish in $\overline{\mathbb{D}}^2$. This became known as the \emph{Shanks Conjecture}.

	The original Shanks Conjecture was disproved shortly after it was formulated. In 1975, Genin and Kamp constructed a two-dimensional least-squares inverse polynomial having a zero in the bidisk \cite{GeninKamp}. Anderson and Jury subsequently proved the conjectured zero-free conclusion for a restricted class of low-degree polynomials, showing that some special cases of the original conjecture do survive \cite{AndersonJury}. Later work simplified and reinterpreted the Genin--Kamp counterexample in the language of modern optimal polynomial approximation \cite{SargentSola}.
	
	However, the polynomial being inverted in the Genin--Kamp example had a zero in the bidisk and it remained plausible that a polynomial that was itself zero-free on the (closed) bidisk would have zero-free OPAs. Accordingly, Delsarte, Genin, and Kamp formulated the following weaker statement in 1980 \cite{DelsarteGeninKampII}.
	
	\begin{conjecture}[Weak Shanks Conjecture]
		\label{conj:weak-shanks}
		Suppose that $f\in\mathbb{C}[z,w]$ has no zeros in $\overline{\mathbb{D}}^2$. Then every optimal polynomial approximant to $1/f$ in $H^2(\mathbb{D}^2)$ is zero-free in $\mathbb{D}^2$.
	\end{conjecture}
	
	We shall call a polynomial $f\in\mathbb{C}[z,w]$ \textit{stable} if
	\[
	f(z,w)\neq 0
	\quad
	\text{for every }(z,w)\in \mathbb{D}^2
	\]
	and \textit{strongly stable} if
	\[
	f(z,w)\neq 0
	\quad
	\text{for every }(z,w)\in\overline{\mathbb{D}}^2.
	\]
	Thus, the Weak Shanks Conjecture asserted that the strong stability of a polynomial forces all of its OPAs to be stable.
	
	The conjecture remained unresolved for more than four decades. Part of its history involved an unsuccessful claimed proof of a modified form of the conjecture by Reddy, Reddy, and Swamy \cite{ReddyReddySwamy}. Delsarte, Genin, and Kamp identified a gap in that argument \cite{DelsarteGeninKampComments}, and the failure was also discussed by Karivaratharajan and Reddy \cite{KarivaratharajanReddy}. Consequently, the weak conjecture remained open.
	
	Several decades later, optimal polynomial approximants were independently developed as a tool for studying cyclicity in Dirichlet-type spaces. B\'en\'eteau, Condori, Liaw, Seco, and Sola studied extremal polynomials and cyclicity in one-variable spaces and on the bidisk \cite{BeneteauCondoriLiawSecoSolaI,BeneteauCondoriLiawSecoSolaII}. 
	Cyclicity of polynomials in Dirichlet-type spaces on the bidisk was characterized by B\'en\'eteau, Knese, Kosi\'nski, Liaw, Seco, and Sola \cite{BeneteauKneseKosinskiLiawSecoSola}; Knese, Kosi\'nski, Ransford, and Sola extended the study to anisotropic Dirichlet-type spaces \cite{KneseKosinskiRansfordSola}.
	B\'en\'eteau, Khavinson, Liaw, Seco, and Sola established strong connections among OPAs, orthogonal polynomials, and reproducing kernels, and proved zero-free results in one-variable Hardy and Dirichlet spaces \cite{BeneteauKhavinsonLiawSecoSola}. B\'en\'eteau, Khavinson, Liaw, Seco, and Simanek subsequently studied the extremal location and asymptotic distribution of zeros of OPAs using Jacobi matrices and Jentzsch-type theorems \cite{BeneteauKhavinsonLiawSecoSimanek}. In particular, their work showed that zeros can enter the unit disk in certain weighted one-variable spaces. This one-variable phenomenon later provided an important model for the resolution of the Weak Shanks Conjecture.
	
	Sargent and Sola developed a general framework for OPAs in multivariable reproducing kernel Hilbert spaces, including their connections with cyclicity, weakly inner functions, orthogonal polynomials, and multidimensional filtering \cite{SargentSola,SargentSolaII}. Their work also highlighted several ways in which the one-variable theory fails to extend directly to higher dimensions. More recently, the second named author established several positive Shanks-type results for restricted classes of two-variable functions, including functions that reduce in an appropriate sense to one-variable problems \cite{FelderShanks}.
	
	The Weak Shanks Conjecture was finally disproved by B\'en\'eteau, Khavinson, and Seco \cite{BeneteauKhavinsonSeco}. They began with the symmetric function, which is zero-free in the closed bidisk, given by
	\[
	F(z,w)=\left(1-\frac{z+w}{\sqrt{6}}\right)^{-5/2}
	\]
	and proved that its affine OPA has a zero in the bidisk. Although $F$ is not a polynomial, it is holomorphic on a neighborhood of the closed bidisk. By approximating $F$ with sufficiently high-degree Taylor polynomials and using continuity of the coefficients of the affine OPA, they obtained a stable polynomial whose affine OPA vanishes in $\mathbb{D}^2$. This produced the first counterexample to the Weak Shanks Conjecture.
	
	The construction of B\'en\'eteau, Khavinson, and Seco was naturally motivated by the earlier one-variable theory. Functions of the form
	\[
	f(z,w)=G\left(\frac{z+w}{2}\right)
	\]
	embed a one-variable weighted problem into the bidisk. On symmetric homogeneous polynomials, the $H^2(\mathbb{D}^2)$-norm induces a one-variable coefficient norm comparable to a weighted Bergman norm. The counterexample was obtained by transferring a one-variable phenomenon, in which zeros of OPAs enter the unit disk in an appropriate weighted space, into the symmetric subspace of the Hardy space on the bidisk \cite{BeneteauKhavinsonLiawSecoSimanek,BeneteauKhavinsonSeco}.
	
	Geronimo and Woerdeman subsequently gave a Fourier-analytic treatment of this phenomenon \cite{GeronimoWoerdeman}. They considered the two-parameter family
	\[
	F_{\alpha,r}(z,w)=\left(1-\frac{z+w}{r}\right)^{-\alpha},
	\qquad r>2,
	\]
	and derived hypergeometric formulas for the Fourier coefficients of $|F_{\alpha,r}|^2$. These formulas permit explicit calculation of the relevant Gram matrix entries and numerical optimization over $\alpha$ and $r$. Their work gives further counterexamples to the Weak Shanks Conjecture and provides a quantitative analysis of the location of the zero variety of the affine OPA. In one optimized example, the affine OPA has a diagonal zero near $(0.975766,0.975766)$.
	Recently, Nav\'io and Seco also constructed symmetric polynomial counterexamples in Dirichlet-type spaces over the bidisk \cite{NavioSeco}.
	
	These counterexamples settle the qualitative question posed by the Weak Shanks Conjecture. They show that the affine OPA associated with a stable polynomial can have a zero variety intersecting the open bidisk. They also lead naturally to a stronger quantitative question: how far into the bidisk can such a zero variety penetrate?
	
	For a polynomial $p\in\mathbb{C}[z,w]$, define
	\[
	\delta_2(p)=\operatorname{dist}_{\mathbb{C}^2}\bigl((0,0),\mathcal{Z}(p)\bigr)
	=\inf\left\{\sqrt{|z|^2+|w|^2}:p(z,w)=0\right\}, 
	\]
	where $\mathcal{Z}(p)$ is the zero variety associated to $p$.
	
	In many one-variable spaces, including $H^2(\D)$, the distance from the origin to any possible zero of a non-constant OPA has a universal positive lower bound; see \cite{BeneteauKhavinsonLiawSecoSola,BeneteauKhavinsonLiawSecoSimanek}. It is natural to ask whether the same holds for $\delta_2(p)$ when $p$ ranges over the affine OPAs associated with stable polynomials. The purpose of the present paper is to show that there is no such lower bound.

	Indeed, our main result is the following:

	\vspace{0.5pc}
	\noindent \textbf{Theorem \ref{thrm:main}.}
	For every $\varepsilon>0$, there exists a stable polynomial $f\in\mathbb{C}[z,w]$ such that, if $p$ is the affine optimal polynomial approximant to $1/f$ in $H^2(\mathbb{D}^2)$, then
	\[
	\operatorname{dist}_{\mathbb{C}^2}\bigl((0,0),\mathcal{Z}(p)\bigr)<\varepsilon.
	\]
	
	The theorem strengthens the existing counterexamples in a significant way; the zero varieties of affine OPAs do not merely enter the bidisk, they may pass arbitrarily close to the origin. In this sense, there is no protected neighborhood of the origin on which the affine OPAs associated with stable polynomials must remain zero-free.
	
	Further, the existence of such functions is given explicitly; we construct a family of polynomials $(f_m)_{m\ge 1}$ for which the affine OPA to $1/f_m$ has a variety passing closest to the origin at the point $(\zeta_m, \zeta_m)$, where
	\begin{equation}
		\zeta_m = \frac{ \|(z+w)f_m\|^2_2  }{ 4\langle f_m, zf_m\rangle}.  
	\end{equation}
	Ultimately, we show that $\zeta_m = O(m^{-1/2})$. Further, knowing the asymptotics, one may deduce that $|\zeta_m| < 1$ for all $m\ge 3$, giving a sequence of counterexamples to the Weak Shanks Conjecture.  
	
	\section{Approach and Outline of Argument}\label{sec:outline}
	The initial counterexamples to the Weak Shanks Conjecture \cite{BeneteauKhavinsonSeco, GeronimoWoerdeman} relied on functions of the form 
	\[
	f(z,w)=G\left(\frac{z+w}{r}\right),
	\qquad r\ge2. 
	\]
	We do not take this approach and our counterexamples are not of this form. Instead, we are motivated by a result proved in \cite{FelderShanks}, which we present here, tailored to the affine OPA setting in the bidisk: 
	\vspace{0.5pc}
	
	\noindent \textbf{Theorem} (\cite[Theorem~7]{FelderShanks}).
	Let $f\in H^2(\D^2)$ be zero-free in the bidisk and let $k_\alpha$ denote the reproducing kernel for $H^2(\D^2)$ at the point $\alpha \in \mathbb{D}^2$. Then the affine OPA to $1/f$ has a zero at $\alpha$ if and only if
	\[
	\det
	\begin{pmatrix}
		\langle k_\alpha, f \rangle & \langle f, zf \rangle & \langle f, wf \rangle \\
		\langle k_\alpha, zf \rangle & \langle zf, zf \rangle & \langle zf, wf \rangle \\
		\langle k_\alpha, wf \rangle & \langle wf, zf \rangle & \langle wf, wf \rangle. 
	\end{pmatrix}
	= 0.
	\]
	\vspace{0.5pc}
	
	Examining the above determinant, one realizes that much simplification is granted by assuming $\langle f, zf \rangle = \langle f, wf \rangle$ and $\langle zf, wf \rangle = \langle wf, zf \rangle$. These equalities can be guaranteed by assuming that $f$ is symmetric (i.e., $f(z,w) = f(w,z)$) and that $f$ has real Fourier coefficients, respectively. 
	The additional constraint of stability makes a counterexample not readily handy. However, realizing that stability of low degree polynomials is invariant under taking integer powers, one is reasonably led to the following family of functions, which serve as our counterexamples. 
	
	For $m\geq 1$, let
	\[
	f_m(z,w) := ((1-a_m z + b_m w)(1+b_mz - a_mw))^m,
	\]
	where
	\[
	r_m = \frac{1}{\sqrt{2m}}, \quad
	a_m = \frac{1-r_m}{2}, \quad
	b_m = \frac{1+r_m}{2}.
	\]
	\vspace{0.5pc}
	Note that each $f_m$ is stable, symmetric, and has real Fourier coefficients. 
	
	Although we do not explicitly verify the determinantal condition above in establishing our main result, it does motivate the properties and form of the functions $f_m$. Our approach begins with the following proposition, the proof of which is given in Section \ref{sec:affine-opa}:
	
	\vspace{0.5pc}
	\noindent \textbf{Proposition \ref{prop: location of closest zeros}.}
	The zero of the affine OPA to $1/f_m$ that lies closest to the origin is $Z_m =(\zeta_m, \zeta_m)$, where
	\begin{equation} \label{eqn:zeta_m}
		\zeta_m = \frac{ \|(z+w)f_m\|^2_2  }{ 4\langle f_m, zf_m\rangle}. 
	\end{equation}
	
	\vspace{0.5pc}
	
	Proving our main result, Theorem \ref{thrm:main}, amounts to showing that $Z_m \to (0,0)$ as $m\to \infty$. Proposition \ref{prop: location of closest zeros} shows that it suffices to establish an upper bound for the numerator of the expression in \eqref{eqn:zeta_m} and a non-zero lower bound for its denominator. These bounds are established in Sections \ref{sec:upper-bound} and \ref{sec:lower-bound}, respectively. We then formally state and prove Theorem \ref{thrm:main} in Section \ref{sec:main-result-proof}.

	\section{The Affine OPA}\label{sec:affine-opa}
	In this section we derive the formula for the affine OPA corresponding to a symmetric function, and then calculate the Euclidean distance between the origin and the variety corresponding to the OPA. 
	
	First, let us observe that if $f \in H^2(\D^2)$ is symmetric, then the affine OPA to $1/f$ must be symmetric. Indeed, if \(f(z,w)=f(w,z)\) and we let $U$ denote the coordinate-swap operator
	\[
	(Ug)(z,w)=g(w,z),
	\]
	then if $p$ is any OPA, we have
	\[
	\|pf-1\|_2 = \|U(pf-1)\|_2 = \|(Up)f-1\|_2.
	\]
	By uniqueness, we have $Up = p$, and thus $p$ must be symmetric.
	
	For any function $f\in H^2(\D^2)$, the optimality of the affine OPA $p$ gives
	\[
	1 - pf \perp \bigvee\left\{f, zf, wf \right\}. 
	\]
	In turn, if the affine OPA is denoted by 
	\[
	p(z, w) = a_0 + a_1z + a_2w,
	\]
	we have 
	\begin{equation}\label{gram}
		\begin{pmatrix}
			\langle f,f\rangle & \langle zf,f\rangle & \langle wf,f\rangle \\
			\langle f,zf\rangle & \langle zf,zf\rangle & \langle wf,zf\rangle \\
			\langle f,wf\rangle & \langle zf,wf\rangle & \langle wf,wf\rangle
		\end{pmatrix}
		\begin{pmatrix}
			a_0\\ a_1 \\ a_2
		\end{pmatrix}
		=
		\begin{pmatrix}
			\langle 1,f\rangle\\ \langle 1,zf\rangle\\ \langle 1,wf\rangle \end{pmatrix}.
	\end{equation}
	Now let us assume that $f$ is symmetric, has real Fourier coefficients, and that $f(0,0)= 1$. Note that in this setting, 
	\[
	\langle f, zf \rangle = \langle f, wf \rangle.
	\]
	and we may denote the affine OPA to $1/f$ by 
	\[
	p(z,w) = A + B(z+w).
	\]

	Simplifying Equation (\ref{gram}) yields
	\[
	\begin{pmatrix}
		\|f\|_2^2 & \langle f,zf\rangle & \langle f,zf\rangle\\
		\langle f,zf\rangle & \|f\|_2^2 & \langle zf,wf\rangle \\
		\langle f,zf\rangle & \langle zf,wf\rangle & \|f\|_2^2 \end{pmatrix} \begin{pmatrix} A\\ B\\ B \end{pmatrix} = \begin{pmatrix} 1\\0\\0 \end{pmatrix}. \]
	Solving this system gives
	\[
	A = \frac{ \|f\|_2^2+\langle zf,wf\rangle }{ \|f\|_2^2\bigl(\|f\|_2^2+\langle zf,wf\rangle\bigr) -2\langle f,zf\rangle^2 }
	\]
	and 
	\[
	B = -\frac{ \langle f,zf\rangle }{ \|f\|_2^2\bigl(\|f\|_2^2+\langle zf,wf\rangle\bigr) -2\langle f,zf\rangle^2 }.
	\]
	We now assume additionally that $\langle f,zf\rangle \neq 0$; otherwise, the affine OPA is constant. 
	
	The zero variety of the OPA is now given by
	\[
	\mathcal Z(p) = \left\{ (z,w)\in\mathbb C^2: z+w=-\frac{A}{B} \right\}.
	\]
	For every \((z,w)\in\mathcal Z(p)\), the Cauchy--Schwarz inequality gives
	\[
	\left|\frac{A}{B}\right|^2 = |z+w|^2 \leq 2\left(|z|^2+|w|^2\right).
	\]
	Consequently,
	\[
	|z|^2+|w|^2 \geq \frac{1}{2}\left|\frac{A}{B}\right|^2.
	\]
	Equality holds if and only if \(z=w\). Since
	\[
	z+w=-\frac{A}{B},
	\]
	the unique point of \(\mathcal Z(p)\) closest to the origin is
	\[
	Z=\left(-\frac{A}{2B},-\frac{A}{2B}\right).
	\]
	Putting $\Delta = \|f\|_2^2 \left( \|f\|_2^2+\langle zf,wf\rangle \right) -2\langle f,zf\rangle^2$ and using
	\[
	A= \frac{\|f\|_2^2+\langle zf,wf\rangle}{\Delta} \qquad\text{and}\qquad B= -\frac{\langle f,zf\rangle}{\Delta},
	\]
	we obtain
	\[
	Z = \left( \frac{\|f\|_2^2+\langle zf,wf\rangle} {2\langle f,zf\rangle}, \frac{\|f\|_2^2+\langle zf,wf\rangle} {2\langle f,zf\rangle} \right).
	\]
	Moreover,
	\[
	\frac12 \|(z+w)f\|_2^2 = \|f\|_2^2+\langle zf,wf\rangle,
	\] and hence
	\[
	Z=(\zeta,\zeta), \qquad \zeta= \frac{\|(z+w)f\|_2^2} {4\langle f,zf\rangle}.
	\]
	In terms of our candidate function $f_m$, we have shown the following: 
	\begin{proposition} \label{prop: location of closest zeros}
		The zero of the affine OPA to $1/f_m$ that lies closest to the origin is $Z_m =(\zeta_m, \zeta_m)$, where
		\[
		\zeta_m = \frac{ \|(z+w)f_m\|^2_2  }{ 4\langle f_m, zf_m\rangle}. 
		\]
	\end{proposition}

	We now compute the asymptotics of $\zeta_m$ by bounding the numerator and denominator in the expression appearing above in Proposition \ref{prop: location of closest zeros}. Ultimately, we obtain that $\zeta_m = O(m^{-1/2})$.

	\section{Upper bound on $\|(z+w)f_m\|^2_2$}\label{sec:upper-bound}
	
	In order to estimate $\|(z+w)f_m\|_2^2$, we split $f_m$ into two factors and bound each factor separately. Recall that
	\[
	f_m(z,w)=\big((1-a_mz+b_mw)(1+b_mz-a_mw)\big)^m,
	\]
	where
	\[
	r_m=\frac{1}{\sqrt{2m}},\qquad
	a_m=\frac{1-r_m}{2},\qquad
	b_m=\frac{1+r_m}{2}.
	\]
	Put
	\[
	g_m(z,w)=(1-a_mz+b_mw)(1+b_mz-a_mw),
	\]
	so that $f_m=g_m^m$. Let $k=\textstyle\lfloor\frac{m}{2}\rfloor$.
	Since $f_m=g_m^kg_m^{m-k}$, we have
	\[
	\|(z+w)f_m\|^2_2
	=\int_{\mathbb{T}^2}|(z+w)g_m^k|^2|g_m^{m-k}|^2\mathop{d\mu}
	\leq
	\|(z+w)g_m^k\|^2_\infty\|g_m^{m-k}\|^2_2.
	\]
	We aim to bound, in terms of $m$, each of the above factors separately. We begin by estimating $\|(z+w)g_m^k\|^2_\infty$.
	
	\begin{proposition}\label{prop for sup numerator bounding}
		There exists a constant $C_1>0$, independent of $m$, such that
		\[
		\|(z+w)g_m^k\|_{\infty}^2\leq C_1\frac{4^k}{m},
		\]
		where $k=\textstyle\lfloor\frac{m}{2}\rfloor$.
	\end{proposition}
	
	\begin{proof}
		Expanding $g_m$ and using $b_m-a_m=r_m$, $a_m+b_m=1$, and
		$a_mb_m=(1-r_m^2)/4$, we obtain
		\[
		g_m(z,w)
		=1+r_m(z+w)+zw-\frac{1-r_m^2}{4}(z+w)^2.
		\]
		By the maximum modulus principle, it suffices to estimate this expression
		on $\mathbb{T}^2$. We shall apply the following elementary observation:
		\[
		zw|z+w|^2
		=zw(z+w)(\overline{z}+\overline{w})
		=(z+w)^2,\quad z,w\in\mathbb{T}.
		\]
		Using this identity, we obtain
		\[
		g_m(z,w)
		=1+r_m(z+w)+zw\left(
		1-\frac{1-r_m^2}{4}|z+w|^2
		\right),\quad z,w\in\mathbb{T}.
		\]
		The factor multiplying $zw$ is nonnegative. Indeed, since
		$|z+w|\leq2$ on $\mathbb{T}^2$, we have
		\[
		1-\frac{1-r_m^2}{4}|z+w|^2
		\geq 1-(1-r_m^2)=r_m^2>0.
		\]
		Applying the triangle inequality, we have
		\begin{align*}
			|g_m(z,w)|
			&\leq 1+r_m|z+w|
			+\left|1-\frac{1-r_m^2}{4}|z+w|^2\right|\\
			&=2+r_m|z+w|-\frac{1-r_m^2}{4}|z+w|^2.
		\end{align*}
		By the elementary inequality
		\[
		\left(\frac{|z+w|}{4}-2r_m\right)^2\geq0,
		\]
		we obtain
		\[
		|g_m(z,w)|
		\leq 2+4r_m^2
		-\left(\frac{3}{16}-\frac{r_m^2}{4}\right)|z+w|^2.
		\]
		By definition of $r_m$,
		\[
		\frac{3}{16}-\frac{r_m^2}{4}\geq\frac{1}{16}.
		\]
		Therefore,
		\begin{equation}\label{eq:gm-torus-bound}
			|g_m(z,w)|
			\leq 2+\frac{2}{m}-\frac{|z+w|^2}{16},
			\quad (z,w)\in\mathbb{T}^2.
		\end{equation}
		At this point, assume $m\geq2$, so that $k\geq1$, and set
		$A=2+\frac{2}{m}$ and $t=|z+w|^2$. The bound
		\eqref{eq:gm-torus-bound} implies that
		\[
		t|g_m|^{2k}\leq t\left(A-\frac{t}{16}\right)^{2k},
		\quad t\in[0,16A].
		\]
		By setting $\Phi(t)=t\left(A-\frac{t}{16}\right)^{2k}$ and using
		elementary calculus, one finds the function $\Phi$ attains its
		maximum at the point $t_0=\frac{16A}{2k+1}$ and decreases afterwards.
		Therefore,
		\[
		\Phi(t_0)
		=\frac{16A}{2k+1}
		\left(\frac{2kA}{2k+1}\right)^{2k}
		\leq \frac{16A^{2k+1}}{2k+1}.
		\]
		Our choice of $k=\textstyle\lfloor\frac{m}{2}\rfloor$ now gives
		$2k\leq m$ and $2k+1\geq m$. Moreover,
		\[
		A^{2k}
		=4^k\left(1+\frac{1}{m}\right)^{2k}
		\leq 4^k e^{2k/m}
		\leq e4^k.
		\]
		In conclusion, we obtain
		\[
		\|(z+w)g_m^k\|_{\infty}^2
		\leq 64e\frac{4^k}{m}.
		\]
	\end{proof}
	
	We now turn toward estimating $\|g_m^{m-k}\|^2_2$. Put $n:=m-k$.
	Using routine algebra, we can write
	\begin{equation}\label{g_m alternate form}
		g_m(z,w)
		=\left(1+r_m\frac{z+w}{2}\right)^2-\frac{(z-w)^2}{4}.
	\end{equation}
	It follows from \eqref{g_m alternate form} that
	\begin{equation}\label{g_m^n expanded}
		g_m^n(z,w)
		=\sum_{j=0}^n\binom{n}{j}
		\left(-\frac{(z-w)^2}{4}\right)^j
		\left(1+r_m\frac{z+w}{2}\right)^{2n-2j}.
	\end{equation}
	We will bound $\|g_m^n\|^2_2$ by splitting the sum
	in \eqref{g_m^n expanded} at $j=\lceil n/4\rceil$, estimating each
	partial sum separately, and applying the triangle inequality. For ease
	of notation, put
	\[
	B_{m,n}(z,w)
	=\sum_{j=0}^{\lceil n/4\rceil-1}\binom{n}{j}
	\left(-\frac{(z-w)^2}{4}\right)^j
	\left(1+r_m\frac{z+w}{2}\right)^{2n-2j}
	\]
	and
	\[
	E_{m,n}(z,w)
	=\sum_{j=\lceil n/4\rceil}^{n}\binom{n}{j}
	\left(-\frac{(z-w)^2}{4}\right)^j
	\left(1+r_m\frac{z+w}{2}\right)^{2n-2j},
	\]
	so that $g_m^n = B_{m,n} + E_{m,n}$. We proceed by providing a sup norm bound on $B_{m,n}$ and a direct Hardy space bound on $E_{m,n}$.

	\begin{lemma} \label{g_m^{m-k} lemma 1}
		There exists a constant $\widetilde{C}>0$, independent of $n = m - k$, such that 
		\[
		\left\|B_{m,n}\right\|_2 \leq \widetilde{C}\frac{2^n}{\sqrt{n+1}}.
		\]
	\end{lemma}
	\begin{proof}
		Using the triangle inequality, and the fact that $\left|z\pm w\right| \leq 1$, we have
		\begin{align*}
			\left\|B_{m,n}\right\|_2 &\leq \sum_{j=0}^{\lceil n/4\rceil-1} \binom{n}{j} \left\|\left( -\frac{(z-w)^2}{4}\right)^j\left( 1+r_m \frac{z+w}{2}\right)^{2n-2j}\right\|_\infty\\
			&\leq \sum_{j=0}^{\lceil n/4\rceil-1} \binom{n}{j} \left( 1+r_m\right)^{2n-2j}\\
			&\leq (1+r_m)^{2n} \sum_{j=0}^{\lceil n/4\rceil-1} \binom{n}{j},
		\end{align*}
		where the last inequality follows from the fact that $1+r_m> 1$. Since $k=\textstyle \lfloor \frac{m}{2}\rfloor$, we have
		\begin{equation} \label{eq:2nr_m^2=1}
			2nr_m^2 = \frac{n}{m} = \frac{m-\lfloor \frac{m}{2}\rfloor}{m} = \frac{ \lceil \frac{m}{2}\rceil}{m} \leq 1
		\end{equation}
		so that $4n^2r_m^2 = 2n\cdot 2nr_m^2 \leq 2n$ and hence $2nr_m \leq \sqrt{2n}$. This inequality, along with the fact that $1+x\leq e^x$ for all $x\in \mathbb{R}$, implies that 
		\begin{equation} \label{g_m^{m-k} lemma 1: eq1}
			\|B_{m,n}\|_2 \leq e^{2nr_m}\sum_{j=0}^{\lceil n/4\rceil-1} \binom{n}{j} = e^{\sqrt{2n}}\sum_{j=0}^{\lceil n/4\rceil-1} \binom{n}{j}.  
		\end{equation}
		We now estimate $\sum_{j<\lceil n/4\rceil} \binom{n}{j}$. Since $3\geq 1$, we have
		\[
		\sum_{j=0}^{\lceil n/4\rceil-1} \binom{n}{j} \leq 3^{n/4}\sum_{j=0}^{\lceil n/4\rceil-1} \binom{n}{j} 3^j.
		\]
		However, the above partial sum is less than the full sum, which is the binomial expansion of $\left( 1+\frac{1}{3}\right)^n$. Therefore,
		\begin{equation}\label{g_m^{m-k} lemma 1: eq2}
			\sum_{j=0}^{\lceil n/4\rceil-1} \binom{n}{j}\leq 3^{ n/4}\sum_{j=0}^n \binom{n}{j} 3^j = 3^{ n/4} \left( 1+\frac{1}{3}\right)^n = \left(\frac{4}{3^{3/4}}\right)^n. 
		\end{equation}
		Combining \eqref{g_m^{m-k} lemma 1: eq1} and \eqref{g_m^{m-k} lemma 1: eq2} gives
		\[
		\|B_{m,n}\|_2 \leq e^{\sqrt{2n}}\left(\frac{4}{3^{3/4}}\right)^n.
		\]
		We claim that the above bound is dominated by $\textstyle \frac{2^n}{\sqrt{n+1}}$. Indeed, 
		\[
		\lim_{n\to \infty} \frac{e^{\sqrt{2n}}\left(\frac{4}{3^{3/4}}\right)^n}{\frac{2^n}{\sqrt{n+1}}} = \lim_{n\to \infty} \sqrt{n+1}e^{\sqrt{2n}}\left( \frac{2}{3^{3/4}}\right)^n=0.
		\]
		Therefore, there exists a constant $\widetilde{C}>0$ (independent of $n$) so that 
		\[
		\|B_{m,n}\|_2 \leq e^{\sqrt{2n}}\left(\frac{4}{3^{3/4}}\right)^n \leq \widetilde{C} \frac{2^n}{\sqrt{n+1}}.
		\]
	\end{proof}
	\begin{lemma} \label{lemma:bound for lj}
		For $j\geq \lceil n/4\rceil$ and $\ell\geq 0$,
		\[
		\left\| \left(\frac{z+w}{2}\right)^\ell \left( - \frac{ (z-w)^2 }{4}\right)^j \right\|^2_2 \leq \frac{\ell!}{(n/2)^\ell}\left\| \left( -\frac{ (z-w)^2 }{4}\right)^j \right\|^2_2.
		\]
	\end{lemma}
	\begin{proof}
		Under the standard coordinate transformation $(z,w)\mapsto (e^{i(\phi + \theta)}, e^{i(\phi-\theta)})$ for $\phi,\theta \in [0,2\pi]$, we can rewrite the norm as
		\[
		\|f\|^2_2 = \frac{1}{4\pi^2} \int_0^{2\pi}\int_0^{2\pi} |f(e^{i(\phi + \theta)}, e^{i(\phi-\theta)})| \mathop{d\phi} \mathop{d\theta}.
		\]
		Under this coordinate transformation,
		\[
		\frac{z+w}{2} = e^{i\phi}\cos(\theta) \hspace{0.25in} \text{and} \hspace{0.25in} -\frac{ (z-w)^2 }{4} = e^{2i\phi}\sin^2(\theta).
		\]
		Therefore, since $|e^{(\ell+2j)i\phi}|=1$,
		\begin{align}
			\left\| \left(\frac{z+w}{2}\right)^\ell \left( - \frac{ (z-w)^2 }{4}\right)^j \right\|^2_2 &= \frac{1}{4\pi^2} \int_0^{2\pi} \int_0^{2\pi} |e^{\ell i\phi}\cos^\ell(\theta)\cdot e^{2ji\phi}\sin^{2j}(\theta) |^2\mathop{d\phi} \mathop{d\theta}\\
			&= \frac{1}{2\pi} \int_0^{2\pi} \cos^{2\ell}(\theta)\sin^{4j}(\theta) \mathop{d\theta}.\label{eq1}
		\end{align}
		Note that, when $\ell=0$, we have
		\begin{equation} \label{eq2}
			\left\| \left( - \frac{ (z-w)^2 }{4}\right)^j \right\|^2_2 = \frac{1}{2\pi} \int_0^{2\pi} \sin^{4j}(\theta) \mathop{d\theta}.
		\end{equation}
		
		Consider the function $F(\theta)=\sin^{4j+1}(\theta)\cos^{2\ell-1}(\theta)$. By the Fundamental Theorem of Calculus, $0=F(2\pi) = \textstyle \int_0^{2\pi} \frac{d}{\mathop{d\theta}}[F(\theta)]\mathop{d\theta}$ and so
		\begin{align*}
			0=(4j+1)\int_0^{2\pi}\sin^{4j}(\theta)\cos^{2\ell}(\theta) \mathop{d\theta} - (2\ell-1)\int_0^{2\pi} \sin^{4j+2}(\theta)\cos^{2\ell-2}(\theta)\mathop{d\theta}.   
		\end{align*}
		Therefore,
		\begin{align*}
			\int_0^{2\pi}\sin^{4j}(\theta)\cos^{2\ell}(\theta) \mathop{d\theta} &=\frac{2\ell-1}{4j+1} \int_0^{2\pi} \sin^{4j+2}(\theta)\cos^{2\ell-2}(\theta)\mathop{d\theta}\\
			&= \frac{2\ell-1}{4j+1}\int_0^{2\pi}\sin^{4j}(\theta)\cos^{2\ell-2}(\theta) (1-\cos^2(\theta))\mathop{d\theta}\\
			&= \frac{2\ell-1}{4j+1}\int_0^{2\pi}\sin^{4j}(\theta)\cos^{2\ell-2}(\theta)\mathop{d\theta}- \frac{2\ell-1}{4j+1}\int_0^{2\pi}\sin^{4j}(\theta)\cos^{2\ell}(\theta)\mathop{d\theta}.
		\end{align*}
		Rearranging gives
		\begin{align*}
			\int_0^{2\pi}\sin^{4j}(\theta)\cos^{2\ell}(\theta) \mathop{d\theta} &= \frac{\frac{2\ell-1}{4j+1} }{1+\frac{2\ell-1}{4j+1}} \int_0^{2\pi}\sin^{4j}(\theta)\cos^{2\ell-2}(\theta)\mathop{d\theta} \\
			&= \frac{2\ell-1}{4j+2\ell}\int_0^{2\pi}\sin^{4j}(\theta)\cos^{2\ell-2}(\theta)\mathop{d\theta}.
		\end{align*}
		Iterating this process gives
		\begin{equation} 
			\int_0^{2\pi}\sin^{4j}(\theta)\cos^{2\ell}(\theta) \mathop{d\theta} = \left(\prod_{i=0}^{\ell-1} \frac{2i+1}{4j+2i+2}\right) \int_0^{2\pi} \sin^{4j}(\theta)\mathop{d\theta}, \label{eq3}
		\end{equation}
		where the $\ell=0$ factor is naturally interpreted to be 1. Observe that, since $j\geq \lceil n/4\rceil \geq n/4$, we have
		\[
		\prod_{i=0}^{\ell-1} \frac{2i+1}{4j+2i+2} \leq  \prod_{i=0}^{\ell-1} \frac{2i+2}{n+2i+2} \leq \prod_{i=0}^{\ell-1} \frac{2i+2}{n} = \frac{2^\ell \cdot \ell!}{n^\ell} = \frac{\ell!}{(n/2)^\ell}.
		\]
		Thus, the formulations given in \eqref{eq1} and \eqref{eq2}, as well as this bound, can be substituted into \eqref{eq3} to obtain the desired result.
	\end{proof}
	
	The following is a well-known combinatorial estimate, the proof of which we include for the sake of completion.
	\begin{lemma} \label{lemma:combinoricsidentity}
		For $q\geq 1$,
		\[
		\frac{4^q}{2\sqrt{q}} \leq \binom{2q}{q} \leq \frac{4^q}{\sqrt{q+1}},
		\]
		where the second inequality also holds for $q=0$.
	\end{lemma}
	\begin{proof}
		We prove the second inequality and note that the first inequality can be shown via a similar argument. One can verify directly that, for $k\geq 0$, we have
		\[
		\left( \frac{k+ \frac{1}{2}}{k+1}\right)^2 \leq \frac{k+1}{k+2}
		\]
		and therefore
		\[
		\frac{ \binom{2k+2}{k+1}  }{  \binom{2k}{k}  } = 4\cdot \frac{k+ \frac{1}{2}}{k+1} \leq 4\sqrt{\frac{k+1}{k+2} }.
		\]
		Thus,
		\[
		\binom{2q}{q} = \prod_{k=0}^{q-1} \frac{ \binom{2k+2}{k+1}  }{  \binom{2k}{k}  } \leq \prod_{k=0}^{q-1} 4\sqrt{\frac{k+1}{k+2} } = \frac{4^q}{\sqrt{q+1}}.
		\]
	\end{proof}
	
	\begin{lemma} \label{lemma:bound on j part}
		For all $n \ge 1$, 
		\[
		\sum_{j=\lceil n/4\rceil}^n \binom{n}{j}^2 \left\| \left( -\frac{ (z-w)^2 }{4}\right)^j \right\|^2_2 \leq \frac{2\cdot 4^n}{n+1}.
		\]
	\end{lemma}
	\begin{proof}
		Using the binomial theorem, one obtains
		\[
		\left\| \left( -\frac{ (z-w)^2 }{4}\right)^j \right\|^2_2 = \frac{1}{16^j} \left\|(z-w)^{2j}\right\|^2_2 = \frac{1}{16^j} \sum_{i=0}^{2j} \binom{2j}{i}^2 = \frac{1}{16^j} \binom{4j}{2j},
		\]
		where the final equality follows from the well-known Vandermonde identity. Since $j\geq \lceil n/4\rceil \geq  n/4$, we have $\textstyle \frac{1}{\sqrt{2j+1}} \leq \frac{\sqrt{2}}{\sqrt{n+2}} \leq \sqrt{\frac{2}{n}}$ and therefore, by putting $q=2j$ into Lemma \ref{lemma:combinoricsidentity}, we further find that 
		\[
		\left\| \left( -\frac{ (z-w)^2 }{4}\right)^j \right\|^2_2 = \frac{1}{16^j} \binom{4j}{2j} \leq \frac{1}{\sqrt{2j+1}} \leq \sqrt{\frac{2}{n}}.
		\]
		Thus, we find that
		\[
		\sum_{j=\lceil n/4\rceil}^n \binom{n}{j}^2 \left\| \left( -\frac{ (z-w)^2 }{4}\right)^j \right\|^2_2 \leq \sqrt{\frac{2}{n}}\sum_{j=\lceil n/4\rceil}^n \binom{n}{j}^2 \leq \sqrt{\frac{2}{n}}\sum_{j=0}^n \binom{n}{j}^2.
		\]
		Another application of Vandermonde's identity and another appeal to Lemma \ref{lemma:combinoricsidentity} guarantees that
		\[
		\sum_{j=0}^n \binom{n}{j}^2=\binom{2n}{n}\leq \frac{ 4^n}{\sqrt{n+1}}
		\]
		so that
		\[
		\sum_{j=\lceil n/4\rceil}^n \binom{n}{j}^2 \left\| \left( -\frac{ (z-w)^2 }{4}\right)^j \right\|^2_2 \leq \sqrt{\frac{2}{n}} \frac{ 4^n}{\sqrt{n+1}} \leq \frac{2\cdot 4^n}{n+1},
		\]
		where the final inequality can be verified directly to hold when $n\geq 1$.
	\end{proof}

	\begin{lemma} \label{g_m^{m-k} lemma 2}
		With $\alpha=\sum_{\ell=0}^\infty \frac{2^\ell}{\sqrt{\ell!}}$,
		\[
		\left\|E_{m,n}\right\|_2 \leq \alpha \sqrt{2}\frac{2^n}{\sqrt{n+1}}.
		\]
	\end{lemma}
	\begin{proof}
		Observe first that, since $\binom{2n-2j}{\ell}=0$ for $\ell>2n-2j$, we can rewrite $E_{m,n}$ in the following way: 
		\begin{align*}
			E_{m,n}(z, w) &= \sum_{j=\lceil n/4\rceil}^n \binom{n}{j} \left( -\frac{(z-w)^2}{4}\right)^j\left( 1+r_m \frac{z+w}{2}\right)^{2n-2j}\\
			&= \sum_{j=\lceil n/4\rceil}^n \left(\binom{n}{j} \left( -\frac{(z-w)^2}{4}\right)^j \sum_{\ell=0}^{2n} \binom{2n-2j}{\ell}r_m^\ell \left( \frac{z+w}{2}\right)^\ell\right)\\
			&= \sum_{\ell=0}^{2n}\left( r_m^\ell \left( \frac{z+w}{2}\right)^\ell\sum_{j=\lceil n/4\rceil}^n \binom{n}{j}\binom{2n-2j}{\ell} \left( -\frac{(z-w)^2}{4}\right)^j \right).
		\end{align*}
		Put
		\[
		P_\ell(z, w) = \left( \frac{z+w}{2}\right)^\ell\sum_{j=\lceil n/4\rceil}^n \binom{n}{j}\binom{2n-2j}{\ell} \left( -\frac{(z-w)^2}{4}\right)^j
		\]
		so that $E_{m,n} = \sum_{\ell=0}^{2n} r_m^\ell P_\ell$. Note that, for a fixed $\ell$, different choices of $j$ yield distinct total degrees $\ell+2j$ so that the summands of $P_\ell$ are orthogonal. Note the following bound:
		\[
		\binom{2n-2j}{\ell} = \frac{ (2n-2j)!}{\ell! (2n-2j-\ell)!} \leq \frac{ (2n-2j)^\ell}{\ell!} \leq \frac{ (2n)^\ell}{\ell!}.
		\]
		We can now get an upper bound on $\|P_\ell\|^2_2$. The orthogonality of the summands of $P_\ell$, as well as Lemmas \ref{lemma:bound for lj} and \ref{lemma:bound on j part} guarantee that
		\begin{align*}
			\|P_\ell\|_2^2 &\leq \sum_{j=\lceil n/4\rceil}^{n}  \binom{n}{j}^2 \binom{2n-2j}{\ell}^2  \left\| \left(\frac{z+w}{2}\right)^\ell \left( - \frac{ (z-w)^2 }{4}\right)^j \right\|^2_2\\ 
			&\leq  \sum_{j=\lceil n/4\rceil}^n \binom{n}{j}^2 \binom{2n-2j}{\ell}^2  \frac{\ell!}{(n/2)^\ell}\left\| \left( -\frac{ (z-w)^2 }{4}\right)^j \right\|^2_2\\
			&\leq  \frac{ (2n)^{2\ell}}{(\ell!)^2}\frac{\ell!}{(n/2)^\ell}  \sum_{j=\lceil n/4\rceil}^n \binom{n}{j}^2  \left\| \left( -\frac{ (z-w)^2 }{4}\right)^j \right\|^2_2\\
			&\leq \frac{2\cdot 4^n}{n+1}  \frac{ (8n)^{\ell}}{\ell!}.
		\end{align*}
		Therefore, $\|P_\ell\|_2 \leq \frac{\sqrt{2}\cdot 2^n}{\sqrt{n+1}} \frac{ \sqrt{8n}^{\ell}}{\sqrt{\ell!}}$. It follows immediately from \eqref{eq:2nr_m^2=1} that $r_m \sqrt{8n}\leq 2$ so that, by taking the full sum,
		\begin{align*}
			\|E_{m,n}\|_2 &\leq \sum_{\ell=0}^{2n} r^{\ell}_m \|P_\ell\|_2\\
			&\leq \frac{\sqrt{2}\cdot 2^n}{\sqrt{n+1}} \sum_{\ell=0}^{\infty} r^{\ell}_m \frac{ \sqrt{8n}^{\ell}}{\sqrt{\ell!}}\\
			&\leq \frac{\sqrt{2}\cdot 2^n}{\sqrt{n+1}} \sum_{\ell=0}^{\infty} \frac{ (r_m\sqrt{8n})^{\ell}}{\sqrt{\ell!}}\\
			&\leq \frac{\sqrt{2}\cdot 2^n}{\sqrt{n+1}} \sum_{\ell=0}^{\infty} \frac{ 2^{\ell}}{\sqrt{\ell!}}.
		\end{align*}
		That $\sum_{\ell=0}^\infty \frac{2^\ell}{\sqrt{\ell!}}$ converges to some finite $\alpha$ follows from a standard Ratio Test argument.
	\end{proof}
	
	\begin{proposition} \label{upper bound for numerator}
		With $k=\lfloor \frac{m}{2} \rfloor$, there exists a constant $C_2 > 0$, independent of $m$, such that
		\[
		\|g_m^{m-k}\|^2_2 \leq C_2 \, \frac{4^{m-k}}{m-k+1}.
		\]
	\end{proposition}
	\begin{proof}
		Recall that $n=m-k$ and use Lemmas \ref{g_m^{m-k} lemma 1} and \ref{g_m^{m-k} lemma 2} to obtain
		\begin{align*}
			\|g_m^n\|_2 &\leq \|B_{m,n}\|_2 + \|E_{m,n}\|_2\\
			&\leq \widetilde{C}\frac{2^{m-k}}{\sqrt{m-k+1}} + \alpha \sqrt{2}\frac{2^{m-k}}{\sqrt{m-k+1}}\\
			&\leq C_2 \frac{4^{m-k}}{\sqrt{m-k+1}}.
		\end{align*}
		That $C_2 = \widetilde{C}+\alpha \sqrt{2}$, the sum of the constants from Lemmas \ref{g_m^{m-k} lemma 1} and \ref{g_m^{m-k} lemma 2}, is independent of $m$ follows from the fact that $\widetilde{C}$ and $\alpha$ are both independent of $m$.
	\end{proof}
	
	All together, we have the following:

	\begin{proposition} \label{final upper bound for numerator} There exists a constant $C>0$, independent of $m$, such that
		\[
		\|(z+w)f_m\|^2_2 \leq C\frac{4^m}{m^2}.
		\]
	\end{proposition}
	\begin{proof}
		Let $k=\textstyle \lfloor \frac{m}{2}\rfloor$. By Propositions \ref{prop for sup numerator bounding} and \ref{upper bound for numerator},
		\begin{align*}
			\|(z+w)f_m\|^2_2 &= \int_{\mathbb{T}^2} |(z+w)g_m^k|^2|g_m^{m-k}|^2 \mathop{d\mu}\\
			&\leq \|(z+w)g_m^k\|^2_\infty \|g_m^{m-k}\|^2_2\\
			&\leq C_1\frac{4^k}{m} \cdot C_2\frac{4^{m-k}}{m-k+1}\\
			&\leq C\frac{4^m}{m^2},
		\end{align*}
		where the last inequality follows from the fact that $m-k+1 \geq m/2$ and setting $C=2C_1\cdot C_2$. 
	\end{proof}
	
	\section{Lower bound on $\langle f_m,zf_m\rangle_{2}$}\label{sec:lower-bound}
	
	We now collect a few technical lemmas necessary for our ensuing estimates:
	
	\begin{lemma} \label{lemma for lower bound}
		For $0\leq j \leq m$, 
		\[
		\left\| \frac{z+w}{2} \left( - \frac{ (z-w)^2 }{4}\right)^j \right\|^2_2 \geq \frac{1}{4m+2} \cdot \frac{1}{2\sqrt{2m}}
		\]
	\end{lemma}
	\begin{proof}
		We use the same coordinate transformation described in the proof of Lemma \ref{g_m^{m-k} lemma 2}. It follows from \eqref{eq3} that
		\begin{align*}
			\left\| \frac{z+w}{2} \left( - \frac{ (z-w)^2 }{4}\right)^j \right\|^2_2 &= \frac{1}{2\pi} \int_0^{2\pi} \cos^2(\theta)\sin^{4j}(\theta)\mathop{d\theta}\\
			&\geq \frac{1}{2\pi} \int_0^{2\pi} \cos^2(\theta)\sin^{4m}(\theta)\mathop{d\theta}\\
			&= \frac{1}{4m+2}\cdot \frac{1}{2\pi} \int_0^{2\pi} \sin^{4m}(\theta)\mathop{d\theta},
		\end{align*}
		where the inequality comes from the fact that $0\leq j\leq m$ and $0\leq \sin^2(\theta)\leq 1$. For $m\geq 1$, a standard integration by parts argument on $\sin^{4m}(\theta)=\sin^{4m-1}(\theta)\sin(\theta)$ and subsequent iteration shows that
		\[
		\frac{1}{2\pi}\int_0^{2\pi} \sin^{4m}(\theta)\mathop{d\theta} = \frac{4m-1}{4m} \frac{1}{2\pi}\int_{0}^{2\pi} \sin^{4m-2}(\theta) \mathop{d\theta} = \prod_{i=1}^{2m} \frac{2i-1}{2i}.
		\]
		All together, we have 
		\begin{align*}
			\left\| \frac{z+w}{2} \left( - \frac{ (z-w)^2 }{4}\right)^j \right\|^2_2 &\geq \frac{1}{4m+2}\cdot \frac{1}{2\pi} \int_0^{2\pi} \sin^{4m}(\theta)\mathop{d\theta}\\
			&=\frac{1}{4m+2}\prod_{i=1}^{2m} \frac{2i-1}{2i}\\
			&= \frac{1}{4m+2}\frac{\binom{4m}{2m}}{16^m}\\
			&\geq  \frac{1}{4m+2}\frac{1}{2\sqrt{2m}}
		\end{align*}
		where the last equality can be verified directly and the last inequality follows from Lemma \ref{lemma:combinoricsidentity}.
	\end{proof}
	
	\begin{lemma} \label{lemma for lower bound 2} For $0\leq j \leq m$,
		\[
		\sum_{j=0}^m (m-j)\binom{m}{j}^2 = \frac{m}{2}\binom{2m}{m}.
		\]
	\end{lemma}
	\begin{proof} 
		Recall that $\binom{m}{j}=\binom{m}{m-j}$ so that $S=\sum_{j=0}^m (m-j)\binom{m}{m-j}^2$. Reindexing via $k=m-j$ gives
		\[
		S = \sum_{j=0}^m (m-j)\binom{m}{m-j}^2 = \sum_{k=0}^m k \binom{m}{k}^2.
		\]
		Thus,
		\[
		2S = \sum_{j=0}^m (m-j)\binom{m}{j}^2 + \sum_{j=0}^m j \binom{m}{j}^2 = (m-j+j)\sum_{j=0}^m  \binom{m}{j}^2 = m\sum_{j=0}^m  \binom{m}{j}^2,
		\]
		from which the desired result follows from dividing by 2 and applying Vandermonde's identity.
	\end{proof}
	
	Finally, we establish our lower bound result on $\langle f_m, zf_m \rangle$:
	
	\begin{proposition} \label{denominator lower bound}
		There exists a constant $D>0$, independent of $m$, such that
		\[
		\langle f_m,zf_m \rangle \geq D\frac{4^m}{m^{3/2}}>0.
		\]
	\end{proposition}
	\begin{proof}
		As noted earlier, it follows from \eqref{g_m alternate form} that
		\begin{equation} \label{g_m^m expanded}
			f_m(z, w) = g_m^m(z,w) = \sum_{j=0}^m \binom{m}{j} \left( -\frac{(z-w)^2}{4}\right)^j\left( 1+r_m \frac{z+w}{2}\right)^{2m-2j}.
		\end{equation}
		Expanding the last factor as we did in the proof of Lemma \ref{g_m^{m-k} lemma 2}, we get
		\[
		f_m(z, w) = \sum_{j=0}^m \sum_{\ell=0}^{2m-2j} \binom{m}{j} \binom{2m-2j}{\ell} r_m^\ell \left( \frac{z+w}{2} \right)^\ell \left( - \frac{ (z-w)^2}{4}\right)^j.
		\]
		Note that $\binom{m}{j} \binom{2m-2j}{\ell} r_m^\ell >0$ in every term. The total degree of a single term of $f_m$ is $2j+\ell$. If we write $\frac{z+w}{2}f_m$ in this way (using separate indices to allow for comparison later), we have
		\[
		\frac{z+w}{2}f_m(z, w) = \sum_{i=0}^m \sum_{k=0}^{2m-2j} \binom{m}{i} \binom{2m-2i}{k} r_m^k \left( \frac{z+w}{2} \right)^{k+1} \left( - \frac{ (z-w)^2}{4}\right)^i,
		\]
		whose individual terms have total degree $2i+k+1$. Consider the inner product $\langle f_m, \frac{z+w}{2}f_m \rangle$. Due to orthogonality, the only non-zero terms of this inner product are those for which the total degrees of the terms of $f_m$ and $\frac{z+w}{2}f_m$ are the same. That is, when $2j+\ell = 2i+k+1$. Put
		\[
		P_{\ell,j}(z, w) = \left( \frac{z+w}{2} \right)^\ell \left( - \frac{ (z-w)^2}{4}\right)^j \quad \text{and} \quad Q_{k,i}(z, w) = \left( \frac{z+w}{2} \right)^{k+1} \left( - \frac{ (z-w)^2}{4}\right)^i.
		\]
		We claim that $\langle P_{\ell,j}, Q_{k,i}\rangle \geq 0$ when $2j+\ell = 2i+k+1$. Using the coordinate transformation described in the proof of Lemma \ref{lemma:bound for lj},
		\begin{equation} \label{eqn: lower bound 1}
			\langle P_{\ell,j} Q_{k,i} \rangle = \frac{1}{2\pi} \int_0^{2\pi} \cos^{\ell+k+1}(\theta)\sin^{2j+2k}(\theta) \mathop{d\theta}.  
		\end{equation}
		Observe that the constraint $2j+\ell = 2i+k+1$ implies that $\ell-k=2(i-j)+1$ is odd so that $\ell+k+1=(\ell-k)+2k+1$, as the sum of two odd integers, is even. Thus, the powers of the integrand in \eqref{eqn: lower bound 1} are even and therefore $\langle P_{\ell,j} Q_{k,i} \rangle \geq 0$. 
		
		All together, we've established that the terms of $\langle f_m, \frac{z+w}{2}f_m \rangle$ are all non-negative. This implies that its value is bounded below by any of its partial sums. We choose the terms arising from $\ell=1$ and $k=0$, where the polynomial factors are the same (i.e., $P_{1,j}=Q_{0,i}$). It now follows from the orthogonality of the $P_{1,j}$ and Lemmas \ref{lemma:combinoricsidentity}, \ref{lemma for lower bound}, and \ref{lemma for lower bound 2} that
		\begin{align*}
			\left\langle f_m, \frac{z+w}{2}f_m \right\rangle &\geq \left\langle \sum_{j=0}^m \binom{m}{j} \binom{2m-2j}{1} r_m P_{1,j}, \sum_{i=0}^m \binom{m}{i} \binom{2m-2i}{0} r_m^0 Q_{0,i} \right\rangle\\
			&= 2r_m \sum_{j=0}^m (m-j)\binom{m}{j}^2 \left\| \frac{z+w}{2} \left( - \frac{ (z-w)^2 }{4}\right)^j \right\|^2_2\\
			&\geq 2r_m \sum_{j=0}^m (m-j)\binom{m}{j}^2 \frac{1}{4m+2}\frac{1}{2\sqrt{2m}}\\
			&= 2r_m \frac{1}{4m+2}\frac{1}{2\sqrt{2m}} \sum_{j=0}^m (m-j)\binom{m}{j}^2\\
			&\geq 2r_m \frac{1}{4m+2}\frac{1}{2\sqrt{2m}} \frac{m}{2}\binom{2m}{m}\\
			&\geq 2r_m \frac{1}{4m+2}\frac{1}{2\sqrt{2m}} \frac{m}{2}\frac{4^m}{2\sqrt{m}}\\
			&=\frac{1}{2m+1} \cdot \frac{4^m}{16\sqrt{m}},
		\end{align*}
		where the last equality follows from the fact that $r_m=\textstyle \frac{1}{\sqrt{2m}}$. One can readily verify that when $f\in H^2(\mathbb{D}^2)$ is symmetric (i.e., $f(z,w)=f(w,z)$), $\langle f, zg\rangle =\langle f, \frac{z+w}{2} f \rangle$. Therefore, 
		\[
		\langle f_m,zf_m \rangle = \left\langle f_m, \frac{z+w}{2}f_m \right\rangle \geq \frac{1}{2m+1} \cdot \frac{4^m}{16\sqrt{m}} \geq \frac{4^m}{48m^{3/2}}>0,
		\]
		where the last inequality uses the basic fact that $2m+1\leq 3m$ for $m\geq 1$. Setting $D=1/48$ concludes the proof.
	\end{proof}
	
	\section{Proof of Main Result}\label{sec:main-result-proof}
	
	\begin{theorem}\label{thrm:main}
		Given any $\epsilon > 0$, there exists a stable polynomial $f\in H^2(\D^2)$ so that
		\[
		\operatorname{dist}_{\C^2}\left(\mathcal{Z}(p), (0,0)\right) < \epsilon,
		\]
		where $\mathcal{Z}(p)$ is the variety associated with the affine OPA to $1/f$.
	\end{theorem}
	\begin{proof}
		For $m\geq 1$, let
		\[
		f_m(z,w) := ((1-a_m z + b_m w)(1+b_mz - a_mw))^m,
		\]
		where
		\[
		r_m = \frac{1}{\sqrt{2m}}, \quad
		a_m = \frac{1-r_m}{2}, \quad
		b_m = \frac{1+r_m}{2}.
		\]
		Each each $f_m$ is stable, symmetric, and has real Fourier coefficients. By Proposition \ref{prop: location of closest zeros}, the zero of the affine OPA to $1/f_m$ that lies closest to the origin is $Z_m =(\zeta_m, \zeta_m)$, where
		\[
		\zeta_m = \frac{ \|(z+w)f_m\|^2_2  }{ 4\langle f_m, zf_m\rangle}. 
		\]
		By Propositions \ref{final upper bound for numerator} and \ref{denominator lower bound}, there exists constants $C, D>0$, independent of $m$, such that
		\[
		\|(z+w)f_m\|^2_2 \leq C\frac{4^m}{m^2} \hspace{0.25in} \text{and} \hspace{0.25in} \langle f_m,zf_m \rangle \geq D\frac{4^m}{m^{3/2}}>0.
		\]
		Therefore,
		\[
		\zeta_m = \frac{ \|(z+w)f_m\|^2_2  }{ 4\langle f_m, zf_m\rangle  } \leq \frac{C}{D} \frac{1}{\sqrt{m}},
		\]
		so that $\zeta_m\to 0$ as $m\to \infty$.
	\end{proof}

	\section{Concluding Remarks}
	We end with a few remarks and observations. First, our counterexamples have the form 
	\[
	f(z,w) = G(z+w, z-w),
	\]
	which may be compared with the form of previously known counterexamples. 
	
	Lastly, a remark of A. Dayan, communicated in \cite{NavioSeco}, suggested that without imposing symmetry, the zero variety of an OPA may approach the origin arbitrarily closely. We note that even while imposing symmetry, this is still true.

	\subsection*{AI Disclosure}
	\small{During the preparation of this manuscript, the authors used OpenAI's ChatGPT (Think Deeper) as an auxiliary tool to explore possible approaches, proof strategies, and mathematical methods relevant to the problems considered in this work, as well as to assist with the refinement of the exposition. All mathematical arguments and proofs were critically examined, independently verified, and refined by the authors. The authors take full responsibility for the mathematical correctness, originality, and integrity of the results and conclusions presented in this manuscript.}
	\subsection*{Acknowledgments} The first named author was financially supported by the National Science Center,
	Poland, under SHENG III research project 2023/48/Q/ST1/00048. Moreover, the first named author acknowledges support for his research stay in University of South Florida from the program Excellence
	Initiative at the Jagiellonian University in Krak\'{o}w, Research Support Module 2026.

	\bibliographystyle{alpha}
	\bibliography{Shanks}

	\Addresses

\end{document}